\documentclass{article}
\usepackage{amssymb}
\begin{document}

\title{\bf A Reconstruction Algorithm for a Semilinear Parabolic Inverse Problem} 

\author{\bf Wuqing Ning, Xue Qin and Yunxia Shang}
\date{}
\maketitle
\begin{center}{\large  School of Mathematical Sciences, University of Science and
Technology of China, Hefei 230026, China }
\end{center}
\begin{center}{Email: wqning@ustc.edu.cn, qx3@mail.ustc.edu.cn, yxshang@mail.ustc.edu.cn}
\end{center}
\space
{\bf Abstract:} {In this paper, we consider an inverse problem to determine a semilinear
 term of a parabolic equation from a single boundary measurement of Neumann type.
 For this problem, a reconstruction algorithm is established by the spectral
 representation for the fundamental solution of heat equation with homogeneous Neumann
  boundary condition and the Whitney extension theorem.}

\section{Introduction}
In this paper, we consider the following initial-boundary value problem
\begin{equation}
\left\{
\begin{array}{ll}
\partial_t u-\Delta u+f(u)=0\ \ &\mbox{in}\ \ \Omega\times(0,T)\\
u(\cdot,0)=0 \ \ &\mbox{in}\ \ \Omega \label{1}\\
u=\varphi \ \ &\mbox{on}\ \ \partial\Omega\times(0,T)\end{array}\right.
\end{equation}
where $\Omega\subset{\Bbb{R}}^n$ $(n\in{\Bbb {N}})$ is a $C^\infty$ bounded
connected domain, $f\in C^1(\Bbb{R})$ is a nondecreasing function satisfying
$f\geq 0,f(0)=0$, and $\varphi\in{C^{2+\alpha,1+\frac{\alpha}{2}}
(\partial\Omega\times[0,T])}$ ($0<\alpha <1$) is such that
$\varphi\geq 0,\varphi(\cdot,0)=0,\varphi\not\equiv 0$.

The semilinear parabolic equation in (\ref{1}) appears, for example, in modeling
enzyme kinetics \cite{ke} and in other models \cite{am,oc}.  Under the above
conditions, it can proved that problem (\ref{1}) has a unique solution $u=u_f$
which is in $C^{2+\alpha,1+\frac{\alpha}{2}}(\overline{\Omega}\times[0,T])$
 \cite{he,la}, where we denote the solution by $u_f(x, t)$ for specifying the
 dependence on the semilinear term $f(\cdot)$. In this paper,  we are concerned
 with the following\\
 \\
{\bf Inverse Problem.} Determine the semilinear term $f=f(\cdot)$ from the
Neumann boundary data $\partial_\nu u_f|_{\partial\Omega\times[0,T]}$.\\
\\
Here $\partial_\nu$ denotes the derivative in the direction of the unit
 outward normal vector to $\partial\Omega$. We should note that,
 for the inverse problem, only boundary data are known and in $\Omega\times(0,T)$
 the solution of problem (\ref{1}) is still unknown.

 Up to now, there have been many research papers on such kind of inverse problem.
For uniqueness results under variant additional assumptions, we refer to
\cite{ca,ch,du,kl,lo}. For existence results, see for instance \cite{pi}. For stability results,
we refer to \cite{ch04,ch06} where the domain is a rectangle in \cite{ch04} and a general smooth
domain is considered in \cite{ch06}. The main ingredients in the proof of the results in
\cite{ch04,ch06} are the maximum principle and the gaussian lower bound for
 the fundamental solution of a certain parabolic operator with Neumann boundary condition.

However, from the point of view of practical applications,
it is more important to give theoretically and numerically reconstruction formula or scheme. To
the authors' knowledge, there is no result about exact reconstruction formula for $f$.
Therefore, the main object of this paper is trying to solve this problem, and the result is
stated as follows.\\
\\
{\bf Theorem 1.} {\em For the inverse problem, it is possible to reconstruct $f$ from
$\varphi$ and $\partial_\nu u_f$ on $\partial\Omega\times[0,T]$.}\\

The proof of Theorem 1 is based on the spectral representation for the fundamental solution of
heat equation with homogeneous Neumann boundary condition and the Whitney extension theorem \cite{wh}
(see e.g. \cite{fe,fe2} for the recent development).
 The reconstruction algorithm established here may be helpful to provide theoretical direction
  for numerical simulation.

\section{Preliminaries}
In order to prove Theorem 1, we show some basic facts as follows.

 First, it is well known that the initial-boundary value problem
\begin{equation}
\left\{
\begin{array}{ll}
\partial_t v-\Delta v=0\ \ &\mbox{in}\ \ \Omega\times(0,T)\\
v(\cdot,0)=0 \ \ &\mbox{in}\ \ \Omega \label{v}\\
v=\varphi \ \ &\mbox{on}\ \ \partial\Omega\times(0,T)\end{array}\right.
\end{equation}
has a unique classical solution $v=v_\varphi(x,t)$.

Now let us recall the definition of the fundamental solution of the heat equation with homogeneous Neumann
 boundary condition. Fix $s\in{(s_0,t_0)}$. Let $U(x,t;y,s)$ be a continuous function in the region:
 $s_0<s<t<t_0$, $x\in\overline{\Omega}$, $y\in\overline{\Omega}$. Then $U(x,t;y,s)$ is called the
 fundamental solution of the heat equation with homogeneous Neumann boundary condition (also called
the heat kernel)
 if, for any
 $u_0(x)\in{C(\overline{\Omega})}$, the function defined by
 $$u(x,t)=\int_\Omega U(x,t;y,s)u_0(y){\rm d}y
 $$
 is a solution of the following initial-boundary value problem
 \begin{equation}
\left\{
\begin{array}{ll}
\partial_t u(x,t)-\Delta u(x,t)=0\ \ &\mbox{in}\ \ \Omega\times(s,t_0)\\ \nonumber
\displaystyle{\lim_{t\rightarrow s}}u(x,t)=u_0(x) \ \ &\mbox{in}\ \ \Omega \\
\partial_\nu u(x,t)=0 \ \ &\mbox{on}\ \ \partial\Omega\times(s,t_0).\end{array}\right.
\end{equation}
In fact, $U(x,t;y,s)$ is the unique solution (see e.g. \cite{it}) of the following
 initial-boundary value problem
\begin{equation}
\left\{
\begin{array}{ll}
\partial_t U(x,t;y,s)-\Delta_x U(x,t;y,s)=0\ \ &\mbox{in}\ \ \Omega\times(s,t_0)\\ \nonumber
\displaystyle{\lim_{t\rightarrow s}}U(x,t;y,s)=\delta(x-y) \ \ &\mbox{in}\ \ \Omega \\
\partial_\nu U(x,t;y,s)=0 \ \ &\mbox{on}\ \ \partial\Omega\times(s,t_0)\end{array}\right.
\end{equation}
where $\delta(\cdot)$ denotes the Dirac delta function.

Next we will give the spectral representation of  the fundamental solution.
 Consider the following eigenvalue problem
\begin{equation}
\left\{
\begin{array}{ll}
-\Delta \omega=\lambda\omega\ \ &\mbox{in}\ \ \Omega\\
\partial_\nu \omega=0\ \ &\mbox{on}\ \ \partial\Omega \label{e}\end{array}\right.
\end{equation}
where $\lambda$ is parameter. As is well known, the collection $\Lambda$ of the
eigenvalues of problem (\ref{e}) is real and countable, and
$\Lambda=\{\lambda_k\}_{k=1}^\infty$ if we repeat each eigenvalue according
 to its finite multiplicity as follows
 $$0=\lambda_1<\lambda_2\leq \lambda_3\leq \lambda_4\cdots\rightarrow +\infty.
 $$
Moreover, there exists an orthonormal basis $\{\omega_k\}_{k=1}^\infty$ of the
 Hilbert space $L^2(\Omega)$, where $\omega_k$ is an eigenfunction corresponding
  to $\lambda_k$:
\begin{equation}
\left\{
\begin{array}{ll}
-\Delta \omega_k=\lambda_k\omega_k\ \ &\mbox{in}\ \ \Omega\\ \nonumber
\partial_\nu \omega_k=0\ \ &\mbox{on}\ \ \partial\Omega. \end{array}\right.
\end{equation}
It can proved that there is a important relationship between the fundamental solution (heat kernel)
  and the spectrum (see e.g. \cite{do}), that is,
\begin{equation}
 U(x,t;y,s)=\sum_{k=1}^\infty e^{-\lambda_k(t-s)}\omega_k(x)\omega_k(y) \label{sp}
\end{equation}
where the convergence is uniform in the region:
 $t-s\geq \varepsilon$, $x\in\overline{\Omega}$, $y\in\overline{\Omega}$ for
  arbitrary $\varepsilon> 0$.

\section{Proof of Theorem 1}
We divide the proof of Theorem 1 into two steps.\\
{\bf Step 1.} First we define an admissible set of unknown semilinear terms by
$$\mathcal{F}=\{f\in C^1(\Bbb{R}): f\geq 0, f(0)=0, f\ \mbox{is nondecreasing}\}.
$$
As is mentioned before, for any $f\in{\mathcal{F}}$, problem (\ref{1}) has a unique solution $u_f\in {C^{2+\alpha,1+\frac{\alpha}{2}}(\overline{\Omega}\times[0,T])}$. Now let $w=w(x,t)=u_f(x,t)-v_\varphi(x,t)$
where $v_\varphi(x,t)$ is the solution of problem (\ref{v}).
 Then by a simple computation we find that $w$ is the unique solution of
 \begin{equation}
\left\{
\begin{array}{ll}
\partial_t w-\Delta w=-f(u_f)\ \ &\mbox{in}\ \ \Omega\times(0,T)\\
w(\cdot,0)=0 \ \ &\mbox{in}\ \ \Omega \label{w}\\
w=0 \ \ &\mbox{on}\ \ \partial\Omega\times(0,T).\end{array}\right.
\end{equation}
In order to express explicitly the solutions of (\ref{v}) and (\ref{w}) in terms of
the fundamental solution $U(x,t;y,s)$, we replace artificially the corresponding boundary conditions by
$\partial_\nu v=\partial_\nu v_\varphi$ and $\partial_\nu w=\partial_\nu u_f-
\partial_\nu v_\varphi$ respectively, which can be done
 because of the solvability of (\ref{1}) and (\ref{v}). Then it follows easily from
 Theorem 9.1 of \cite{it} (take $\alpha=\beta=0$) that
 $v_\varphi$ and $w$ can be expressed by
\begin{equation}
v_\varphi(x,t)=\int_0^t\int_{\partial\Omega} U(x,t;y,s)
 \partial_\nu v_\varphi(y,s){\rm d}S(y){\rm d}s  \label{ev}
\end{equation}
 and
\begin{equation}
\begin{array}{ll}
&w(x,t)=\displaystyle{-\int_0^t\int_\Omega} U(x,t;y,s)f(u_f(y,s)){\rm d}y{\rm d}s\\
&\ \ \ +\displaystyle{\int_0^t\int_{\partial\Omega}} U(x,t;y,s)[\partial_\nu u_f-
\partial_\nu v_\varphi](y,s){\rm d}S(y){\rm d}s  \label{ew}
\end{array}
\end{equation}
for $x\in{\overline{\Omega}}$ and $t\in{[0,T]}$, where ${\rm d}S(y)$ is the surface
 element on $\partial\Omega$. Here we have made use of  the fact that
   $U(x,t;y,s)$ satisfies the homogeneous Neumann boundary condition in
   $(y,s)$ (see \cite{it} pp.59).  Noting that $v=\varphi$ and $w=0$
on $\partial\Omega\times[0,T]$, we obtain by (\ref{ev}) and (\ref{ew}) that,
for $(x,t)\in{\partial\Omega\times[0,T]}$,
\begin{equation}
\begin{array}{ll}
I&:=\displaystyle{\int_0^t\int_\Omega} U(x,t;y,s)f(u_f(y,s)){\rm d}y{\rm d}s\\
&=\displaystyle{\int_0^t\int_{\partial\Omega}} U(x,t;y,s)[\partial_\nu u_f-
\partial_\nu v_\varphi](y,s){\rm d}S(y){\rm d}s\\
&=\displaystyle{\int_0^t\int_{\partial\Omega}} U(x,t;y,s)\partial_\nu u_f
(y,s){\rm d}S(y){\rm d}s-\varphi(x,t)\\
&:=a(x,t).  \label{a}
\end{array}
\end{equation}
Since $\varphi$ and $\partial_\nu u_f$ on $\partial\Omega\times[0,T]$ are already known,
so is $a(x,t)$.\\
{\bf Step 2.}  In this step we reconstruct $f(\cdot)$. For fixed $s\in[0,T]$, from Section 2
we may expand $f(u_f(\cdot,s)$ in $L^2(\Omega)$ by
\begin{equation}
f(u_f(y,s))=\sum_{k=1}^\infty c_k(s)\omega_k(y).  \label{f}
\end{equation}
Consequently, in view of (\ref{sp}), for $x\in\overline{\Omega}$, we have
\begin{equation}
I=\sum_{k=1}^\infty p_k(t)\omega_k(x), \ \ p_k(t)=\int_0^t  e^{-\lambda_k(t-s)}c_k(s){\rm d}s. \label{p}
\end{equation}
On the other hand, by the standard Whitney extension theorem (see e.g. \cite{wh}), for fixed
$t\in[0,T]$ there exists an extension
${\tilde{a}}(\cdot,t)$ of the function $a(\cdot,t)$ from $\partial\Omega$ to ${\Bbb{R}}^n$.
Set $a_k(t)=({\tilde{a}}(\cdot,t),\omega_k(\cdot))_{L^2(\Omega)}$ where $(\cdot,\cdot)_{L^2(\Omega)}$
denotes the inner product in $L^2(\Omega)$. In view of (\ref{a}) and (\ref{p}), we may expand formally $a$
in the same way. Therefore, if we
 let $p_k(t)=a_k(t)$, then by (\ref{p})
it is easy to see that
$c_k(\cdot)=(a'_k(\cdot)+\lambda_k a_k(\cdot))$. Hence by (\ref{f}) we have
\begin{equation}
f(u_f(y,s))=\sum_{k=1}^\infty (a'_k(s)+\lambda_k a_k(s))\omega_k(y).  \label{f2}
\end{equation}
For any $x\in{\partial\Omega}$, there exists a sequence $\{y_m\}_{m=1}^\infty$ such that
$\Omega\ni y_m\rightarrow x$ as $m\rightarrow\infty$. Passing to the limit and Noting $u_f=\varphi$
on $\partial\Omega\times[0,T]$,
we obtain by (\ref{f2})
 \begin{equation}
f(\varphi(x,s))=\sum_{k=1}^\infty (a'_k(s)+\lambda_k a_k(s))\omega_k(x)\ \ \mbox{on}\ \ \partial\Omega\times[0,T].  \label{f3}
\end{equation}
The right hand side of (\ref{f3}) is determined by $a$ and hence by $\varphi$ and $\partial_v u_f$ on
$\partial\Omega\times[0,T]$. And if we let $x$ and $s$ vary on $\partial\Omega\times[0,T]$, then the values
of $\varphi$ vary on the interval
$[0,\max_{\partial\Omega\times[0,T]}\varphi]$. Thus by using the known data we can establish a map from
$[0,\max_{\partial\Omega\times[0,T]}\varphi]$ to $\Bbb{R}$, which gives the desired function $f(\cdot)$. $\hfill\square$\\
\\
{\bf Remark.} We point out here that the function $f$ reconstructed above is independent of the extension ways
of the function $a$. In fact, this assertion can be proved by the previous results on uniqueness and
 stability for the inverse problem (see Section 1).

\section{Conclusion}
In this paper, we have established a reconstruction algorithm for an inverse problem to determine a semilinear
 term of a parabolic equation from an additional single boundary measurement. The key idea is to find an
  intrinsic relation between the known data and the unknown function. Although the algorithm  formula can not
   be explicitly written out, in theory it is still important for solving practical problem. The strategy used
   here may be effective for other similar semilinear inverse problems or even for more challenging nonlinear
   inverse problems (see e.g. \cite{is,ka,ni}).
\section*{Acknowledgement}

This work was partially supported by the National Natural Science Foundation of China (Grant No. 11101390), the fundamental Research Funds for the Central University and JSPS Fellowship P05297.

\end{document}